# The asymptotic representation of some series and the Riemann hypothesis

By

## M. Aslam Chaudhry[*] and Gabor Korvin[**]

*Department of Mathematics and Statistics
**Department of Earth Sciences
King Fahd University of Petroleum and Minerals
Dhahran 31261, Saudi Arabia.
E-mail: maslam@kfupm.edu.sa, gabor@kfupm.edu.sa

**Abstract**

We present a conjecture about the asymptotic representation of certain series. The conjecture implies the Riemann hypothesis and it would also indicate the simplicity of the non-trivial zeros of the zeta-function

## 1. Introduction

We present a conjecture about the asymptotic behavior of a class of functions represented by a series. A special case of the conjecture has applications in the theory of zeta functions ([1, 2, 3, 4, 5]).

Riemann [6] proved that the zeta-function

$$\zeta(s) := \sum_{n=1}^{\infty} \frac{1}{n^s} \qquad (s = \sigma + i\tau, \sigma > 1), \qquad (1.1)$$

has a meromorphic continuation to the complex plane which satisfies the functional equation ([5], p.13, eq. (2.1.1))

$$\zeta(s) = 2(2\pi)^{-(1-s)} \cos(\frac{\pi}{2}(1-s))\Gamma(1-s)\zeta(1-s) := \chi(s)\zeta(1-s). \qquad (1.2)$$

The Riemann zeta function has simple zeros at $s = -2, -4, -6, \ldots$ called *trivial zeros*. All other zeros, called *non-trivial zeros*, of the zeta function are symmetric about the *critical line* $\sigma = 1/2$ and are in the *critical strip* $0 \leq \sigma \leq 1$. The multiplicity of these non-trivial zeros (in general) is not known. Riemann conjectured that all non-trivial zeros of the zeta function lie on the critical line $\sigma = 1/2$. This conjecture is called the *Riemann Hypothesis* (RH, for short). We present a conjecture on the asymptotic behavior of a class of functions. It will be shown that the proof of a special case of the conjecture would imply the truth of RH and demonstrate the simplicity of the zeros of the zeta function. For definitions and terminology, we refer to [7, 8, 9].

## 2. A restricted result about the asymptotic representation of series

We start with a trivial observation. Let $u(x)$ ($x \geq x_0 > 0$) be a decreasing, continuous and positive function such that $\sum_{n=1}^{\infty} (-1)^{n-1} u(x+n)$ ($x \geq x_0 > 0$) is uniformly convergent.

Then we note that

$$s_2(x) := u(x+2) - u(x+1) \geq 0, \tag{2.1}$$

and

$$s_{2n}(x) := \sum_{k=1}^{2n} (-1)^k u(x+k) = s_{2n-2}(x) + (u(x+2n-1) - u(x+2n)) \geq s_{2n-2}(x) \geq 0, \tag{2.2}$$

and

$$\begin{aligned} s_{2n}(x) = u(x+1) - [(u(x+2) - u(x+3)) + (u(x+4) - u(x+5)) \\ + \ldots + (u(x+2n-2) - u(x+2n-1)) + u(x+2n)] \leq u(x+1) \leq u(x). \end{aligned} \tag{2.3}$$

Therefore, the sequence $\{s_{2n}(x)\}$ of the even partial sums is increasing and bounded. Hence, the sequence is convergent and satisfies the inequality

$$\sum_{n=1}^{\infty} (-1)^{n-1} u(x+n) = \lim_{n \to \infty} s_{2n}(x) \leq u(x) \qquad (\forall x > 0), \tag{2.4}$$

which shows that

$$\sum_{n=1}^{\infty}(-1)^{n-1}u(x+n)=O(u(x)) \qquad (x\to\infty). \qquad (2.5)$$

One does not find a similar result about the asymptotic behavior of a function defined by the series $\sum_{n=1}^{\infty}u(x+n)$. It seems very unlikely to have a general result about the asymptotic of such series. However, we have a restricted result: Let us define $f(x) \sim g(x)$ $(x\to\infty)$ if $\dfrac{f(x)}{g(x)} \to 1$ $(x\to\infty)$.

**Lemma**

Let $u(x)$ $(x\geq x_0 >0)$ be a continuous function such that $u(x) \sim \dfrac{1}{x^s}$ ($s=\sigma+i\tau$, $\sigma>0$, $x\to\infty$) and $\sum_{n=1}^{\infty}u(x+n)$ $(x\geq x_0 >0)$ is uniformly convergent. Then

$$\sum_{n=1}^{\infty}u(x+n)=O(x^{1-\sigma}) \qquad (x\to\infty). \qquad (2.6)$$

**Proof** Assume that $0<\sigma\leq 1$. Since the series is uniformly convergent this implies the existence of an integer $N$ such that $\sup_x\{\sum_{n>N}u(x+n)\}<1$, and $\sup_x\{\sum_{n=1}^{N}u(x+n)\}<C$. Therefore, the series $\sum_{n=1}^{\infty}u(x+n)$ is bounded and trivially of the size $O(x^{1-\sigma})$. Secondly, assume that $\sigma>1$. In this case, we have $\left|\sum_{n=1}^{\infty}u(x+n)\right| \leq C_1\sum_{n=1}^{\infty}(x+n)^{-\sigma} \leq C_2 x^{1-\sigma}$ and we are done.

We give an example where the asymptotic relation (2.6) can be verified:

Consider $\sum_{n=1}^{\infty}u(x+n)$ ($x\geq x_0 >0$), where

$$u(x) := x^{-s} \qquad (\sigma > 1). \qquad (2.7)$$

In this particular case,

$$\sum_{n=1}^{\infty} u(x+n) = \frac{1}{\Gamma(s)} \int_0^{\infty} \frac{t^{s-1}}{e^t - 1} e^{-xt} dt. \qquad (2.8)$$

However,

$$\frac{1}{\Gamma(s)} \frac{t^{s-1}}{e^t - 1} \sim \frac{1}{\Gamma(s)} t^{s-2} \qquad (t \to 0^+). \qquad (2.9)$$

Therefore, by Watson's lemma ([8], p. 5),

$$\sum_{n=1}^{\infty} u(x+n) = \frac{1}{\Gamma(s)} \int_0^{\infty} \frac{t^{s-1}}{e^t - 1} e^{-xt} dt \sim \frac{\Gamma(s-1)x^{-s}}{\Gamma(s)} = \frac{x^{1-s}}{s-1} \qquad (x \to \infty), \qquad (2.10)$$

what is in agreement with the (2.6).

## 3. The Main conjecture and its motivation

Consider the function $u(x) := x^{-s}$ ($s := \sigma + i\tau$, $\sigma > 1$), let $\mu(n)$ denote the Möbius function ([4], p. 217), and let $(s)_0 := 1$, $(s)_k := \Gamma(s+k)/\Gamma(s)$ be the Pochhammer symbol. Then

$$\sum_{n=1}^{\infty} \mu(n) u(x+n) = \sum_{n=1}^{\infty} \frac{\mu(n)}{(n+x)^s} = \sum_{k=0}^{\infty} \frac{(s)_k (-x)^k}{k!} \left(\sum_{n=1}^{\infty} \frac{\mu(n)}{n^{k+s}}\right) = \sum_{k=0}^{\infty} \frac{(s)_k (-x)^k}{\zeta(s+k) k!}$$

$$(0 \le x < 1, \sigma > 1), \qquad (3.1)$$

which extends the well-known identity ([4], p. 260(1))

$$\sum_{n=1}^{\infty} \frac{\mu(n)}{n^s} = \frac{1}{\zeta(s)} \qquad (\sigma > 1). \qquad (3.2)$$

The RHS of (3.1) can be represented as an inverse Mellin transform ([4], p. 224):

$$\sum_{n=1}^{\infty} \mu(n) u(x+n) = \sum_{k=0}^{\infty} \frac{(s)_k (-x)^k}{\zeta(s+k) k!} = \frac{1}{2\pi i} \int_{c-i\infty}^{c+i\infty} \frac{\Gamma(z)\Gamma(s-z)}{\Gamma(s)\zeta(s-z)} x^{-z} dz$$

$$(\sigma >1,\ 0<c_1 \leq c \leq c_2 <\sigma -1). \tag{3.3}$$

Using the asymptotic relation ([1], p. 6, eq. (1.45))

$$\Gamma(s) = \Gamma(\sigma \pm i\tau) = \sqrt{2\pi}\,|\tau|^{\sigma -1/2}\exp(-\frac{\pi}{2}|\tau|)(1+O(1/|\tau|))$$

$$(|\tau|\to\infty, -\infty < a \leq \sigma \leq b < \infty), \tag{3.4}$$

for the gamma function, we have

$$|\Gamma(z)\Gamma(s-z)| \sim (2\pi)|y|^{\sigma -1}e^{-\pi |y|} \qquad (|y|\to\infty). \tag{3.5}$$

In the integrand of (3.3), as we have $\text{Re}(s-z)\geq \sigma_0 >1$ and the fact that the zeta function does not vanish in the region $\sigma >1$, the factor $1/\zeta(s-z)$ remains bounded in the region $0<c_1 \leq c \leq c_2 <\sigma -1$. Consequently, in view of the asymptotic representation (3.5) and due to the fact that the factor $1/\zeta(s-z)$ remains bounded in the region $0<c_1 \leq c \leq c_2 <\sigma -1$, the inverse Mellin transform integral (3.3) exists, and it can be evaluated by Cauchy's residue theorem. *The presence of the exponential term in (3.5) shows that the error term tends to zero faster than any power as* $y\to\infty$ ([9], p.148). The poles $z=-n$ ($n=0,1,2,...$) of the integrand in (3.3) are on the LHS of the line of integration. The integrand also has poles to the RHS of the line of integration at $z=s+n$ ($n=0,1,2,...$) and at $z=s-\rho$, where the $\rho$'s are the non-trivial zeros of the zeta function. Taking the sum over the residues gives, for large values of $x$ ([9], p. 148),

$$\sum_{n=1}^{\infty} \mu(n)u(x+n) \sim \sum_{n=0}^{\infty} R_n x^{-s-n} + \sum_{\rho} \mathfrak{R}_\rho x^{\rho -s} \qquad (x\to\infty), \tag{3.6}$$

where

$$R_n := x^{s+n}\,\text{Res}[\frac{\Gamma(z)\Gamma(s-z)}{\Gamma(s)\zeta(s-z)}x^{-z}; s+n], \tag{3.7}$$

$$\mathfrak{R}_\rho := x^{s-\rho} \operatorname{Re} s\left[\frac{\Gamma(z)\Gamma(s-z)}{\Gamma(s)\zeta(s-z)} x^{-z}; s-\rho\right]. \tag{3.8}$$

It is to be remarked that *if all zeros of the zeta function are simple then $\mathfrak{R}_\rho$ does not depend on x and* hence the asymptotic behavior of (3.6) is governed by the power of *x with the largest positive real part*, that is,

$$\sum_{n=1}^{\infty} \mu(n) u(x+n) = O(x^{\sigma_M - \sigma}) \qquad (x \to \infty), \tag{3.9}$$

where $\sigma_M := \sup\{\operatorname{Re}\rho : \zeta(\rho) = 0\}$.

**Conjecture**  Let $u(x)$ $(x \geq x_0 > 0)$ be a continuous function such that $u(x) \sim \dfrac{1}{x^s}$

$(s := \sigma + i\tau, \ \sigma > 1, \ x \to \infty)$ and $\sum_{n=1}^{\infty} u(x+n)$ $(x \geq x_0 > 0)$ is uniformly convergent. Then

$$\sum_{n=1}^{\infty} \mu(n) u(x+n) = O(x^{\frac{1}{2}-\sigma}) \qquad (x \to \infty), \tag{3.10}$$

where $\mu(n)$ is the Möbius function.

This conjecture would imply RH and resolve the simplicity of the zeros of the zeta function. It is because, for the particular case $s = 3/2 + i\tau$, we would have

$$\sum_{n=1}^{\infty} \frac{\mu(n)}{(n+x)^{\frac{3}{2}+i\tau}} = O\left(\frac{1}{x}\right) \qquad (x \to \infty), \tag{3.11}$$

that is,

$$\sigma_M := \sup\{\operatorname{Re}\rho : \zeta(\rho) = 0\} = 1/2. \tag{3.12}$$

It is to be remarked that the particular case (3.11) of the conjecture, once proven, would prove the Riemann hypothesis and resolve the simplicity of the zeros of the zeta function.

## 4. Discussion and concluding remarks

There are well over one hundred known equivalents or consequences of RH. (Recent reviews [10, 11] list some fifty, covering a span of a century, from Helge von Koch's classic 1901 paper [12] on the connection between RH and the *Prime Number Theorem* to Baez-Duarte's new versions of the *Nyman-Beurling Criterion* in 2002 [13].) Our conjecture is stronger than RH; its particular case (for $\sigma = \frac{3}{2}$) is equivalent to RH. Because of the simple formulation and analytic tractability of this conjecture, we envisage that it would serve as a launching pad for many new attacks on RH. A less drastic conjecture may be stated under the conditions of the conjecture as:

$$\sum_{n=1}^{\infty} \mu(n) u(x+n) = O(x^{\frac{1}{2}-\sigma+\varepsilon}) \qquad (x \to \infty, \ \forall \varepsilon > 0). \qquad (4.1)$$

**Acknowledgements**   The authors are grateful to the King Fahd University of Petroleum and Minerals for the excellent research facilities.